\documentclass[11pt]{article}
\usepackage{amsmath, amssymb}

\newtheorem{thm}{Theorem}

\newtheorem{lm}[thm]{Lemma}

\renewcommand{\phi}{\varphi}
\renewcommand{\epsilon}{\varepsilon}

\newcommand{\BB}{\mathbb}
\newcommand{\g}{\mathfrak}
\newcommand{\pf}{\noindent {\it Proof. }}
\newcommand{\qed}{$\qquad$ $\square$ \vskip7pt}
\newcommand{\separate}{\vskip7pt}

\begin{document}

\title{The $K$-group of $\BB R$-constructible Sheaves}
\author{Matvei Libine}
\maketitle

\begin{abstract}
Let $Z$ be a smooth projective manifold.
In these notes I will prove that the $K$-group of $\BB R$-constructible
sheaves is isomorphic to the free abelian group with one generator for
each open semialgebraic subset $U$ (which I will denote by the same
letter) modulo the Mayer-Vietoris relations:
$U+V - U \cup V - U \cap V =0$.

I will prove it by showing that both groups in question are isomorphic
to the group of all integer-valued semialgebraic functions on $Z$.
\end{abstract}

Let us introduce notations first. Let $Z$ be a smooth projective manifold.
The $K$-group of $\BB R$-constructible sheaves on $Z$
will be denoted by $\g K$.
Let $\g G$ denote the factor-group of
the free abelian group with one generator for each open semialgebraic
set $U$ by the subgroup generated by expressions of the form
$U+V - U \cup V - U \cap V =0$. Let $\g F$ denote the group of all
integer-valued semialgebraic functions on $Z$.
For any subset $S \subset Z$, let $\chi_S$ denote the indicator
function of $S$, $\BB C_S$ the constant sheaf on $S$, and
$j_S: S \hookrightarrow Z$ the inclusion map.

\begin{thm}
$\g G \simeq \g F \simeq \g K$. The composition isomorphism
$\g G \tilde \to \g K$ sends $U$ into $[j_{U!}(\BB C_U)]$.
\end{thm}

\pf
Since for any open $U, V \subset Z$ we have
$\chi_U + \chi_V - \chi_{U \cup V} - \chi_{U \cap V} =0$,
we have a well-defined map
$$
\alpha: \g G \to \g F, \qquad U \mapsto \chi_U.
$$
Integer-valued semialgebraic functions are integral linear
combinations of characteristic functions of semialgebraic
sets. By definition, a semialgebraic set is a disjoint
union of locally closed ones. Let $Q \subset Z$ be a
locally closed semialgebraic set, write $Q=U \cap C$, where
$U$ and $C$ are open and closed semialgebraic sets respectively.
Then
$$
\chi_Q=\chi_{U \cap C}=\chi_U - \chi_{U \setminus C}
$$
which shows that $\alpha$ is surjective.

\begin{lm} $\alpha$ is injective and hence isomorphism.
\end{lm}

\pf
Suppose that $\sum_i a_i U_i \in G$ ($a_i \in \BB Z\setminus\{0\}$)
is such that $\sum_i a_i \chi_{U_i}=0$. We need to show that
$\sum_i a_i U_i=0$ in $\g G$.
Applying the Mayer-Vietoris relations we can assume that
all the $U_i$'s are connected. Grouping terms we can assume
that all the $U_i$'s are different. Let $M=\max_i |a_i|$.

First we will show that using the Mayer-Vietoris relations
we can arrange the $U_i$'s so that
any two $U_j$, $U_k$ are either disjoint or one is contained
in the other. Let us call two sets $A$ and $B$ \emph{overlapping}
if $A \cap B \ne \varnothing$, $A \setminus B \ne \varnothing$
and $B \setminus A \ne \varnothing$. That is $A$ and $B$ overlap
if they have nonempty intersection and none of them is contained
in the other.

Fix a finite Borel measure $\mu$ on $Z$ such that $\mu(U)>0$
whenever $U$ is a nonempty open set.

Let us check that we can always find two overlapping sets with
coefficients of the same sign as long as there are some
overlapping sets.
Suppose that $U_j$ and $U_k$ overlap and $a_j>0$, $a_k<0$.
Pick points $x_1 \in U_j \cap \partial U_k$ and
$x_2 \in U_j \cap \partial U_k$ (such points exist because $U_j$
and $U_k$ are connected). Then
$(a_j\chi_{U_j}+a_k\chi_{U_k})(x_1)>(a_j\chi_{U_j}+a_k\chi_{U_k})(x_2)$.
Since $\sum_i a_i \chi_{U_i}(x_1)=\sum_i a_i \chi_{U_i}(x_2)=0$,
there exists $l$ such that $a_l \chi_{U_l}(x_1) < a_l \chi_{U_l}(x_2)$.
Let us say $a_l>0$, then $x_2 \in U_l$ and $x_1 \notin U_l$.
It easily follows that $U_j$ and $U_l$ overlap and $a_j$ has same sign
as $a_l$. If $a_l<0$, then for the same reasoning $U_k$ and $U_l$
overlap and $a_k$, $a_l$ have same signs.

Suppose that  $U_j$ and $U_k$ overlap and $a_j \ge a_k >0$
or $a_j \le a_k <0$, then we replace $a_jU_j + a_kU_k$ with
$(a_j-a_k)U_j+a_k (U_j \cup U_k) + a_j (U_j \cap U_k)$.
The open set $U_j \cap U_k$ may not be connected in which case
we replace it with algebraic sum of its connected components.
Observe that the number $\sum_i |a_i|\mu(U_i)$ does not increase,
and it decreases if and only if at least one cancellation is
performed. That is either $U_j \cup U_k$ or one of the connected components
of $U_j \cap U_k$ appeared in the original expression $\sum_i a_i U_i$
with coefficient of opposite sign.
Also, the number $\max_i |a_i|$ does not increase either.

We will show that if we keep performing these replacements then the
process will stop after a finite number of steps, that is we will
arrive at an expression with no overlaps.

Indeed, observe that the number of expressions of the form
$\sum_i b_i V_i$ where $0<|b_i| \le M$ and $V_i$'s are distinct
connected open sets obtained from $U_i$'s by taking unions and
connected components of intersections is finite.
Thus to show that the process stops it is enough to check that
we cannot get the same expression twice.
For the sake of contradiction suppose an expression $E$
appeared twice. Then because
the function $\sum_i |a_i|\mu(U_i)$ is non-increasing it must
stay constant. Thus there will be no cancellations in the
coefficients at each step.

Let $\Omega$ consists of all open
sets $U$ satisfying the following two properties:
\begin{enumerate}
\item The set $U$ appears in the expression $E$ or in an intermediate
expression.
\item If $U$ appears in $E$, then there is an intermediate
expression where $U$ appears with a coefficient of larger
absolute value than in $E$.
\end{enumerate}
$\Omega$ is not empty because the very first step introduces an element
of $\Omega$, namely a connected component of $U_j \cap U_k$.
Clearly, $\Omega$ is finite. Let $V$ be a minimal element of $\Omega$,
that is $V$ is not contained in any other $V' \in \Omega$.
But then once $V$ appears in an intermediate expression with
coefficient of greater absolute value than in $E$ it cannot
be canceled. Thus we arrive at contradiction which proves
that no expression can appear twice. Therefore, the process
terminates at an expression with no overlaps.

Our last step is to show that the sum $\sum_i a_i U_i$ is empty
when there are no overlaps. Suppose not, then pick any $U_j$ which
is not contained in any other $U_k$.
Since $a_j \ne 0$ there must be $U_k$ such that
$U_k \varsubsetneq U_j$ and no other $U_l$ contains $U_k$.
Since $U_j$ is connected, $U_j \cap \partial U_k \ne \varnothing$
and we can pick a point $x \in U_j \cap \partial U_k$.
Then because there are no overlaps
$0=\sum_i a_i \chi_{U_i}(x)= a_j \chi_{U_j}(x)=a_j \ne 0$,
a contradiction. This finishes our proof of the lemma.
\qed

\separate

Next we construct a map
$$
\beta: \g K \to \g F.
$$
If $\S$ is an $\BB R$-constructible sheaf, $\beta(\S)$ is a function
whose value at a point $z \in Z$ is the dimension of the stalk of $\S$
at $z$. Clearly, this map descends to the $K$-group.
It is also clear that if $U$ is an open semialgebraic subset of $Z$,
then $\beta(j_{U!}(\BB C_U))= \chi_U$.
On the other hand, the elements $[j_{U!}(\BB C_U)]$ generate $\g K$
and one can show that the map $\beta$ is injective and surjective
using the same reasoning as for $\alpha$. Thus $\beta$ is an isomorphism.
Alternatively, Theorem 9.7.1 in \cite{KaScha} states that $\beta$ is
an isomorphism.

\separate

As a result we obtain an isomorphism
$$
\beta^{-1} \circ \alpha : \g G \to \g K
$$
which sends $U$ into $[j_{U!}(\BB C_U)]$.
\qed

\noindent
{matvei@math.umass.edu}

\noindent
{DEPARTMENT OF MATHEMATICS AND STATISTICS, UNIVERSITY OF MASSACHUSETTS,
LEDERLE GRADUATE RESEARCH TOWER, 710 NORTH PLEASANT STREET, AMHERST,
MA 01003}

\enddocument